\documentclass[reqno]{amsart} 
\usepackage{latexsym}
\usepackage{amsmath}
\usepackage{amsthm}
\usepackage{graphicx}
\begin{document} 

\title[\hfilneg \hfil Fels sculptures]
{The mathematics of Fels sculptures} 

\author[ D. Fels and A. B. Mingarelli]
{David Fels and Angelo B. Mingarelli}  
\address{Artist in Residence, School of Architecture,\\ Carleton University, Ottawa, Ontario, Canada, K1S 5B6}
\email[D. Fels]{dfelsm@gmail.com}
\address{School of Mathematics and Statistics\\ 
Carleton University, Ottawa, Ontario, Canada, K1S\, 5B6}
\email[A. B. Mingarelli]{amingare@math.carleton.ca}

\date{February 2, 2014}
\thanks{Submitted February, 2014.}
\subjclass[2000]{00A71, 00A69}
\keywords{Sculpture, Fels sculptures, Applications of mathematics}

\begin{abstract}
We give a purely mathematical interpretation and construction of sculptures rendered by one of the authors, known herein as Fels sculptures. We also show that the mathematical framework underlying Ferguson's sculpture, {\it The Ariadne Torus}, may be considered a special case of the more general constructions presented here. More general discussions are also presented about the creation of such sculptures whether they be virtual or in higher dimensional space.
\end{abstract}

\maketitle
\section*{Introduction}
Sculptors manifest ideas as material objects. We use a system wherein an idea is symbolized as a solid, governed by a set of rules, such that the sculpture is the expressed material result of applying rules to symbols. The underlying set of rules, being mathematical in nature,  may thus lead to enormous abstraction and although sculptures are generally thought of as three dimensional objects, they can be created in four and higher dimensional (unseen) spaces with various projections leading to new and pleasant three dimensional sculptures. 

The interplay of mathematics and the arts has, of course, a very long and old history and we cannot begin to elaborate on this matter here. Recently however, the problem of creating mathematical programs for the construction of ribbed sculptures by Charles Perry was considered in \cite{hs}. For an insightful paper on topological tori leading to abstract art see also \cite{chs}. Spiral and twirling sculptures were analysed and constructed in \cite{ea}. This paper deals with a technique for sculpting works mostly based on wood (but not necessarily restricted to it) using abstract ideas based on twirls and tori, though again, not limited to them.

\begin{figure}[ht!]
\centering
\includegraphics[width=53mm]{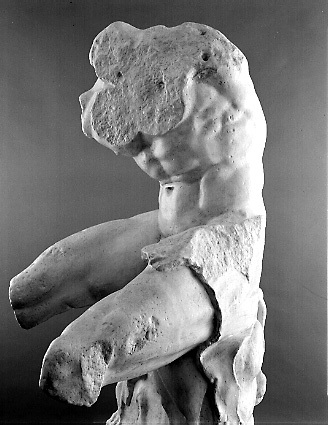}
\includegraphics[width=36mm]{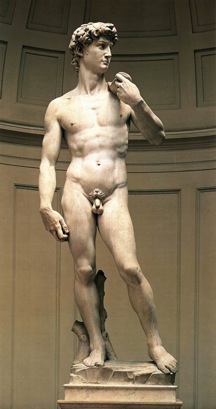}
\includegraphics[width=50mm]{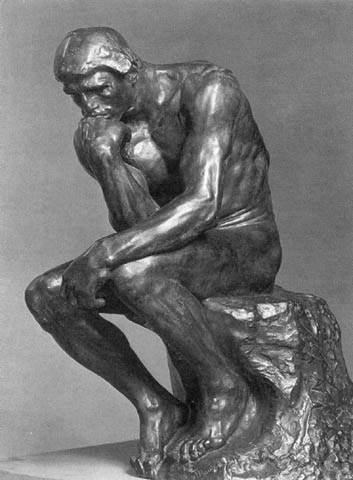}
\caption{{\it Apollononios} son of Nestor, Belvedere Torso, 1st century BC, INV. 1192, Vatican Museum, Rome; {\it David} by Michelangelo, 1501-04, Galleria dell'Accademia, Florence, Italy; August Rodin, {\it The Thinker}, 1902 Bronze, Mus\'ee Rodin, Paris.
}
\label{Fels6}
\end{figure}

The construction is based on the idea that a solid toroidally equivalent structure can be built in such a way that its cross-sections are curvilinear quadrangles that twirl, rotate and dilate as one follows a typical arm from its starting point and then back to it. The more basic Fels scuptures (see Figure 1, left) are reminiscent of the work of Helaman Ferguson \cite{hf} where the cross sections are now curvilinear triangles (the {\it Ariadne torus})  though the surfaces are usually encoded with a surface filling curve. Such sculptures are the work of the first author who has been creating these for over 25 years. A sculptor has a choice of using an existing language of form or creating one’s own specific set of rules. The ancient Greeks, Michelangelo, and Rodin (Figure~\ref{Fels6}) all employed the relationships inherent in the human body as a pre-existing set of rules when defining their language of sculptural forms. 

\begin{figure}[ht!]
\centering
\includegraphics[width=70mm]{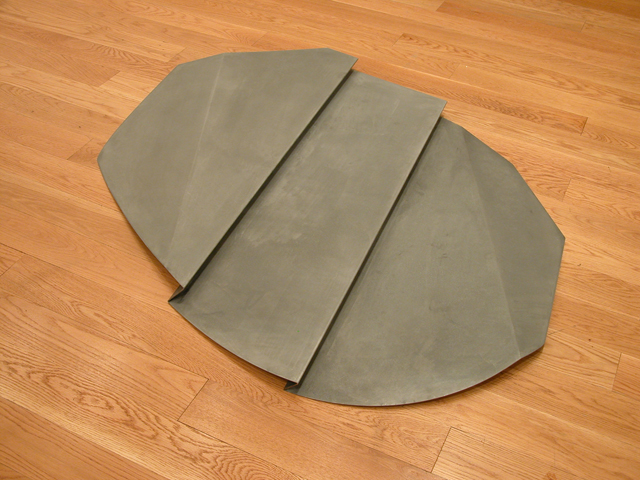}
\caption{{\it Phantom Group: Open Conic Plane, 2 Parallel Double Systems, 2 Parallel Single Systems}, David Rabinovitch, 1967.}
\label{Fels7}
\end{figure}

David Rabinovitch (Figure~\ref{Fels7}) took the rules of sculptural language to the extreme by creating sculptures simply by sending rules to a fabricator. For example, a steel fabricator would be given the rule,``Make it higher on the right” as the complete language of form for the fabrication of a sculpture. 

\begin{figure}[ht!]
\centering
\includegraphics[width=50mm]{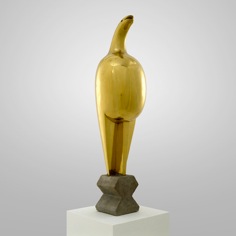}
\includegraphics[width=33.5mm]{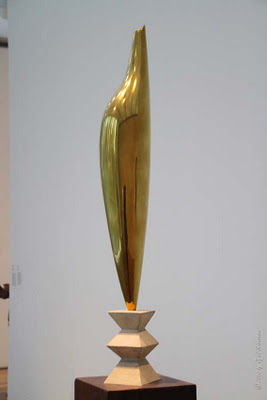}
\includegraphics[width=50mm]{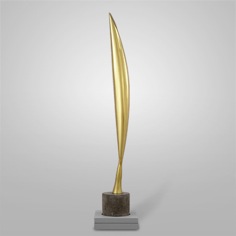}
\caption{Constantin Brancusi, {\it Maiastra}, 1912? Guggenheim Museum, New York; {\it Golden Bird}, 1919-20, Art Institute of Chicago; {\it Bird in Space}, 1932, Guggenheim Museum, New York.}
\label{Fels9}
\end{figure}

In the mid 20th century Brancusi (Figure~\ref{Fels9}) used reductionism to finding a personal language of form; wherein multiple versions of the same sculpture gradually revealed the simplest language of form for the expression of Brancusi’s concept or idea.

\begin{figure}[ht!]
\centering
\includegraphics[width=50mm]{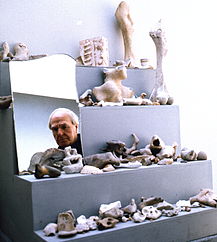}\quad\quad\quad
\includegraphics[width=42mm]{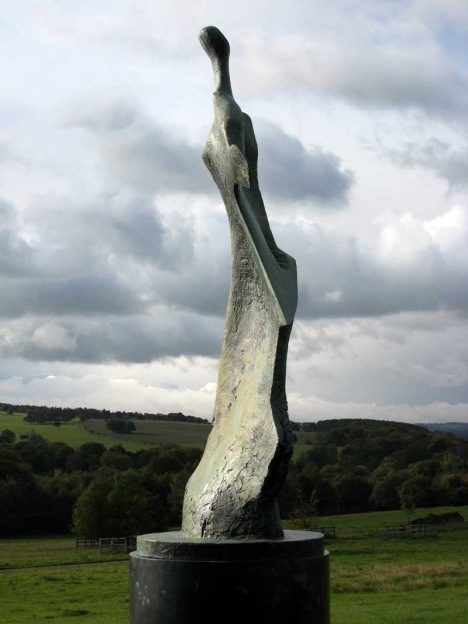}
\caption{{\it Henry Moore in his studio} with a collection of naturally 
found objects, England, 1975, photo by Allan Warren; {\it Standing Figure Knife Edge}, Bronze, 1961.}
\label{Fels10}
\end{figure}

Henry Moore derived the rules of his sculptural forms from careful examination of natural objects in conjunction with the human form, resulting in understandable abstractions (Figure~\ref{Fels10}).

\begin{figure}[ht!]
\centering
\includegraphics[width=100mm]{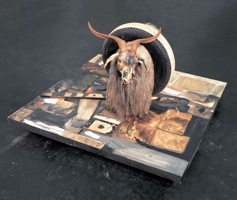}
\caption{Robert Rauschenberg, Monogram, 1955-1959, free standing combine, Moderna Museet, Stockholm}
\label{Fels8}
\end{figure}

More recently sculptors have been creating works of such multifarious materials and combinations of images and ideas that the language or languages of form employed in the work are complex and hard to identify (e.g., Figure~\ref{Fels8})

In his sculptures Fels searched and found a language of form balanced between innate human feelings and understanding of three-dimensional form and a rigorous simple set of underlying rules; wherein the rules are consistent from sculpture to sculpture but do not limit the expression of complex ideas. For example the sculpture ‘Sailing Through Time’ was constructed using the mathematical ideas elucidated below. 

\begin{figure}[ht!]
\centering
\includegraphics[width=90mm]{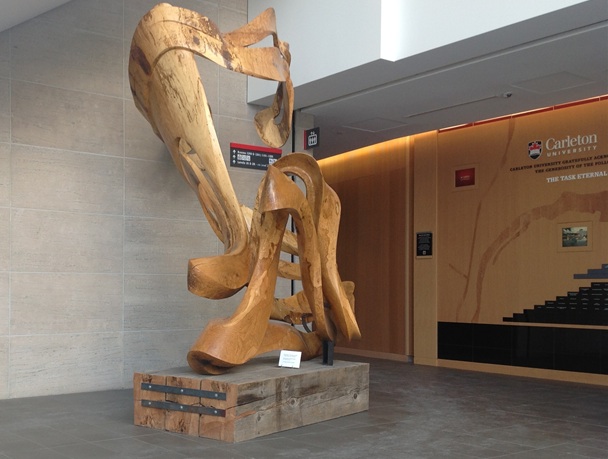}
\caption{{\it Sailing Through Time} by David Fels, 2012 (courtesy of Carleton University, Canada)}
\label{Fels4}
\end{figure}

Fels' rules pay homage to human experience, for example the twists of the sculptures relate to movement and change in human attitude, captured in a set of rules amenable to mathematical scrutiny unique to both fields. Moreover, as discussed in this paper, Fels' language of sculptural form mathematically leads to a dramatic expansion of the concept of sculpture in higher dimensions and directly enables real-time interactive sculptural expression and creation on a virtual pallet.

\section*{Curvilinear quadrangles}
Let $A, B, C, D$ be four distinct points in the euclidean plane and assume that each pair of points $AB$, $BC$, $CD$ and $DA$ is  joined by a piecewise differentiable {\it non-self-intersecting}\footnote{this is not, strictly speaking, adhered to in this paper but represents the simplest of cases.} curve (ususally convex or concave in this paper). The resulting planar figure is then called a {\it curvilinear quadrilateral} or {\it curvilinear quadrangle} and denoted by the symbol, $\mathcal{R}$. It follows that a curvilinear quadrilateral is a closed Jordan curve \cite{pm} having a finite area.

Since, in the sequel, all such quadrangles are composed of natural materials (e.g., wood) curvilinear quadrilaterals will have a mass density function $\rho(x,y)$ that can be used to define its {\it center of mass},  $(\overline{x}, \overline{y})$,  where the coordinates $\overline{x}$, $\overline{y}$ are given by the usual
$$\overline{x} = \frac{1}{m}\iint\limits_ {\mathcal R} x\rho(x,y)\, dA,\quad \overline{y} = \frac{1}{m}\iint\limits_ {\mathcal R} y\rho(x,y)\, dA$$
where $m$, the mass of the quadrangle, is given by 
$$m = \iint\limits_ {\mathcal R} \rho(x,y)\, dA$$
and $dA$ is an element of area in euclidean coordinates, \cite{tan}.

In the absence of a mass-density function we replace the center of mass concept by the point intesection of the main diagonals of the quadrangle under investigation.

\section*{The construction of the simplest Fels sculptures}
Now, let $\mathcal{C}$ denote a closed non-self-intersecting differentiable curve in euclidean space, $\mathbf{R}^3$, parametrized by the (position) vector field ${\mathbf{r}}(t) = x(t)\mathbf{i} + y(t) \mathbf{j}+z(t)\mathbf{k}$, where, without loss of generality, we may assume that $t$ lies in some finite interval $[0,a] := 0 \leq t \leq a$ (where $a>0$ is a fixed positive real number). Since $\mathcal{C}$ is closed, it is necessary that ${\mathbf{r}}(0)={\mathbf{r}}(a)$.

We define a sequence of {\it planar} curvilinear quadrangles $\mathcal{Q}(t)$, for $t$ in $[0,a]$ with the property that the center of mass of  $\mathcal{Q}(t)$ coincides with ${\mathbf{r}}(t)$, for a given parameter value, $t$. (This property is an essential feature of a {\it Fels sculpture}.) It is also necessary that all four vertices of such a typical quadrangle be coplanar. We will assume, in addition, that the sequence or family, $\mathcal{Q}(t)$, of such quadrangles deforms continuously (in the sense of topology) for our $t$. Intuitively, this means that as one moves from a quadrangle $\mathcal{Q}(t_1)$ to $\mathcal{Q}(t_2)$ along the curve $\mathcal{C}$ one can expect that their areas, perimeters, etc. be close provided that $t_1$ and $t_2$ are close. 

There is one more caveat, however. Said deformation is also assumed to {\it shrink} (i.e., areas decrease) and {\it rotate} a quadrangle $\mathcal{Q}(t)$, for $t>0$, in an initially clockwise (or initially counterclockwise) direction until, {\it in the simplest of cases}, a single  distinguished value  $t=t^*$  is reached, where $0 < t^* < a$ (corresponding to a maximum rotation by an angle $\pi/\alpha$, where $0<\alpha<\infty$, and $\alpha$ is determined uniquely by $t^*$) and $\mathcal{Q}(t^*)$ has minimal area amongst all quadrangles in the family. When $t^*<t<a$ the deformation proceeds in the same initial direction so that the areas of the $\mathcal{Q}(t)$ now increase and, in the end, $\mathcal{Q}(0)= \mathcal{Q}(a)$ exactly (i.e., the vertices all match as do the quadrangles). The resulting, necessarily, non-self-intersecting  solid structure $\mathcal{F}$ in three-dimensional space, consisting of the union of all such quadrangles $\mathcal{Q}(t)$ for $t$ in the interval $[0,a]$, is called a (simple) {\it Fels sculpture}.

\begin{figure}[ht!]
\centering
\includegraphics[width=50mm]{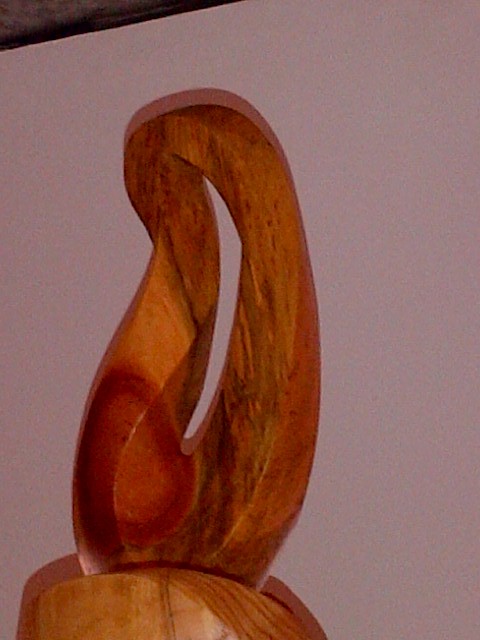}
\includegraphics[width=90mm]{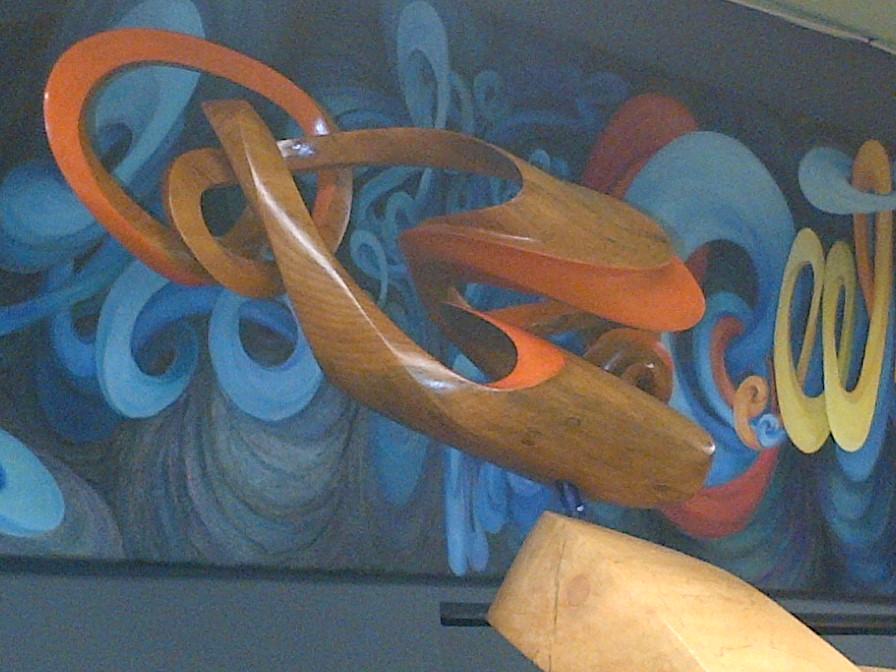}
\caption{Simple Fels scultpures}
\label{Fels1}
\end{figure}

Thus, roughly speaking, a simple Fels sculpture can be thought of as a three-dimensional toroidal structure (topologically equivalent to a doughnut or torus) with surface cross sections being curvilinear quadrangles of decreasing and increasing areas, amongst which there is only one cross section having minimal area (and thus only one having maximal area, by construction). In other words, for a given simple Fels sculpture there is a unique cross-section having minimal area and a unique cross-section having maximal area. (This latter observation becomes a {\it requirement} below when considering complex Fels sculptures.)

\begin{figure}[ht!]
\centering
\includegraphics[width=50mm]{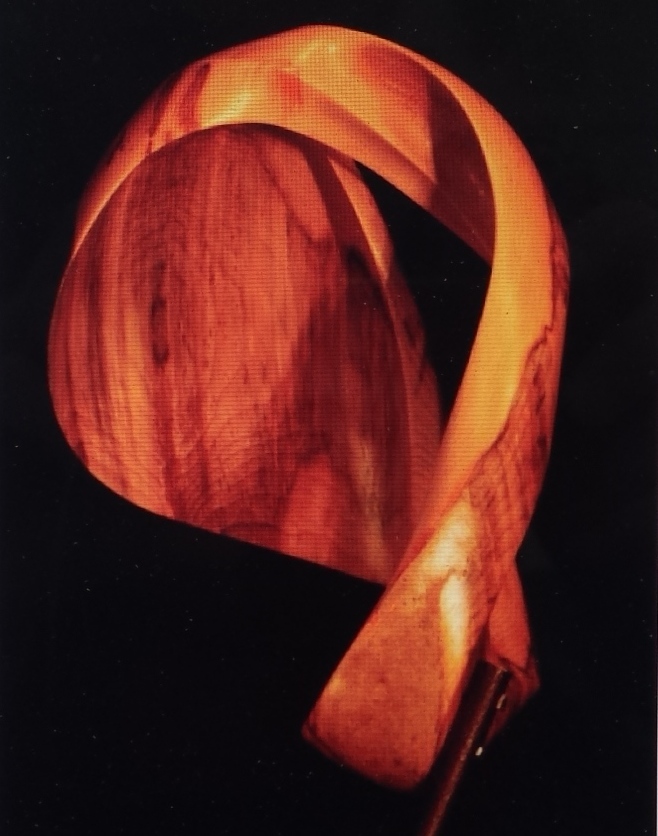}
\caption{Another Fels sculpture}
\label{Fels1}
\end{figure}

\begin{figure}[ht!]
\centering
\includegraphics[width=50mm]{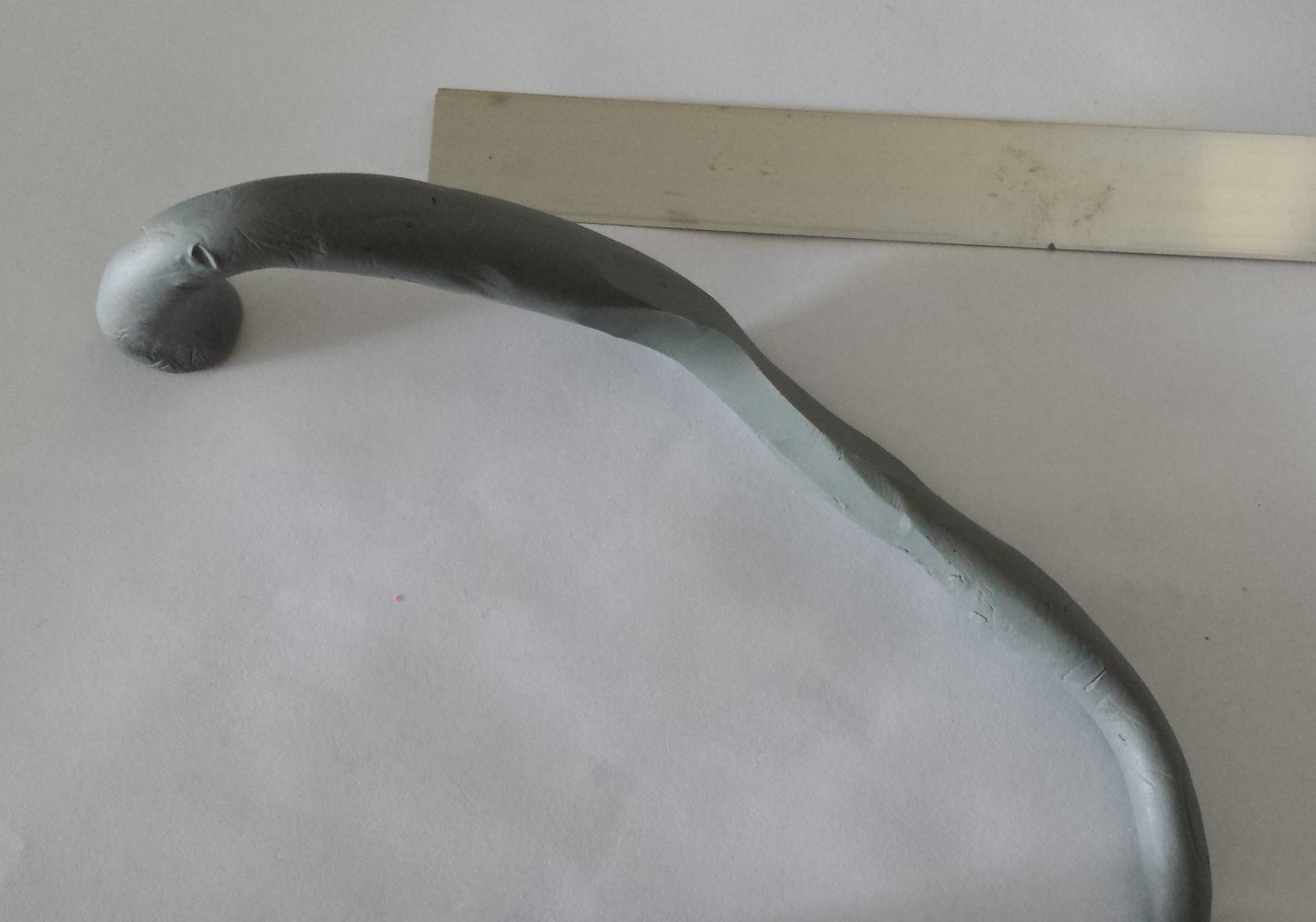}
\caption{Tubular frame and quadrangular shavings as the basis for a Fels sculpture}
\label{Fels18}
\end{figure}

In practice there is much insight to be gained by starting with circles (instead of quarangles) whose centers lie on the curve $\mathcal{C}$ and satisfy the same conditions, as figures, as do the quadrangles. This will give a tubular solid. One can then {\it shave off} four sets of arcs to create quadrangles of the desired shape thus completing the Fels sculpture (Figure~\ref{Fels18}).

Of course, the above construction applies to three-dimensional simple Fels sculptures though there is nothing exceptional in the underlying space being of dimension equal to three. The construction presented here, being purely mathematical, may be extended to any finite number of dimensions to form Fels sculptures in hyperspace (see the end). Depending on the angle being viewed these can then be projected back into three dimensions to give a new (angle-dependent from four dimensions) Fels sculpture in three dimensional euclidean space.

{\bf EXAMPLE}

\begin{figure}[ht!]
\centering
\includegraphics[width=90mm]{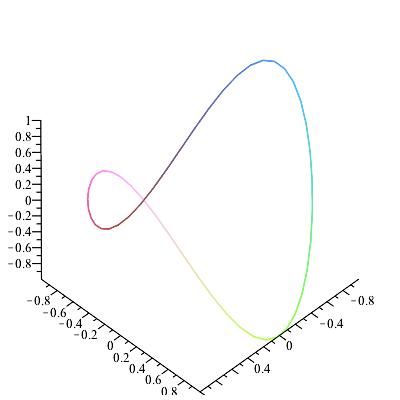}
\caption{The closed differentiable curve $\mathbf{r}(t) = (\cos t, \sin t, \cos 2t ), \  0 \leq t \leq 2\pi$, from which a Fels scultpure may be derived. }
\label{Fels12}
\end{figure}

We now show specifically how such a simple Fels sculpture can be created (virtually or {\it on paper}) from first principles. Let's start with the three dimensional curve, $\mathcal{C}$, (Figure~\ref{Fels12}) given parametrically on the interval $[0,a] = [0,2\pi]$ by defining its position vector to be 
$$\mathbf{r}(t) = (\cos t, \sin t, \cos 2t ), \quad 0 \leq t \leq 2\pi.$$
This is a closed non-self-intersecting differentiable curve in $\mathbf{R}^3$ with the property that for any fixed $t$, the normal vector, $\mathbf{n}$, at the point $P(\cos t, \sin t, \cos 2t )$ is given by $$\mathbf{n} = (\cos t, 5\sin t, 1).$$ 
The plane cross-section to $\mathcal{C}$ at a point $t=t_0$ is (by definition) given by the plane that contains both the point $P$ and the straight line, $\mathcal{L}$, through $P$ parallel to $\mathbf{n}$. In this case, we see that $\mathcal{L}$ may be described parametrically by the equations (in cartesian coordinates)
\begin{eqnarray*}
x &=& \cos t_0 + s \cos t_0,\\
y & = & \sin t_0 + 5s\sin t_0\\
z&=&\cos 2t_0+s,
\end{eqnarray*}
where the parameter $s$ varies over the real numbers. Since the tangent vector, $\mathbf{r}^{\prime}(t_0)$, to $\mathcal{C}$ and the normal to the plane $\Pi$ coincide, the equation of this plane cross section $\Pi$ through $P$ at $t=t_0$ is given by
$$\Pi:\quad  (\sin t_0)\, x - (\cos t_0)\, y + (2\sin 2t_0)\, z = \sin 4t_0.$$
Next, the angle, $\alpha$, $0 \leq \alpha \leq \pi/2$, between this plane $\Pi$ and the plane $z=0$ (i.e., the $xy$-plane) is given by
$$\alpha = \frac{\pi}{2} - \theta,$$
where $\cos \theta = \mathbf{n}\cdot \mathbf{i}/|\mathbf{n}|,$ and $\mathbf{i}$ is the usual unit vector in the $x$-direction. Substituting our value for $\mathbf{n}$ we get
$$\alpha = \frac{\pi}{2} - \arccos\left (\frac{\sqrt{2}\cos t_0}{2\sqrt{1+12\sin^2t_0}}  \right).$$
Now the euclidean sphere 
$$ (x-\cos t_0)^2 + (y-\sin t_0)^2 + (z-\cos 2t_0)^2 = R^2$$
centered at $P$ intersects the plane $\Pi$ in a circle of radius $R$ whose projection onto the plane $z=0$ is given by the elliptical curve
$$16(x-\cos t_0)^2 + 16(y-\sin t_0)^2 + ((\csc t_0)y - (\sec t_0)x)^2 = R^2.$$
A continuous deformation of the type required by the above construction may now be seen by letting $R$ diminish with increasing $t$. The {\it speed} of the deformation is unimportant here, just that such a deformation exists, and this is now clear. For example, we can let $R$ vary with $t$ according to, say, $R = (1+t)^{-1}$ up to a critical value  $t=t^*$ where we now produce a deformation that forces $R$ to {\it increase} for $t^*< t < 2\pi$ until the final/original circle is reached. (Note that rotations of the basic circular figure here lead nowhere due to the complete symmetry of the circle with respect to rotations.) The resulting solid is a closed tubular surface.

\section*{Ferguson's Ariadne Torus}

We show briefly how the construction of Ferguson's {\it  Ariadne torus} \cite{hf}, denoted briefly by AT,  without the incised space filling curve, can essentially be subsumed in the present context. Although we do not have an original Ferguson sculpture before us to compare angles, etc. the idea should be clear. 

\begin{figure}[ht!]
\centering
\includegraphics[width=90mm]{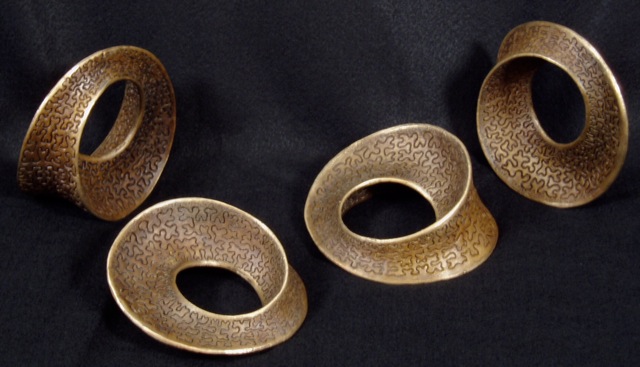}
\caption{The Ariadne Torus, by Helaman Ferguson.}
\label{Fels2}
\end{figure}

Consider a vertical cross section $\Pi$ of the AT (defined as the plane $\Pi$ containing the normal vector to a point $P$ on the boundary of AT and orthogonal to the tangent vector at $P$), cf., Figure~\ref{Fels2}. The intersection of $\Pi$ with the AT defines a hyperbolic triangle, $\mathcal{T}$, maybe even an {\it equilateral} hyperbolic triangle (although this is a guess). The curve containing all the centers of mass of these cross-sections is clearly closed and appears to have an elliptical projection onto a plane. The triangle $\mathcal{T}$ is now dilated and rotated up to a maximum angle of roughly $2\pi/3$, thus making $\alpha=3/2$ in the above discussion on Fels sculptures. Now choose an $a>0$. Since the final ``triangle" does not match the original one (as it is a rotation of the original $\mathcal{T}$) this meets all the requirements of a simple Fels sculpture except for the assumption that the final triangle match the original triangle exactly (no rotations of the initial figure are allowed in the Fels case).

\section*{Curvilinear polygons and complex Fels sculptures}

We can imagine that there is no need for the restriction that the base cross-sections of our three-dimensional sculptures be either triangles (as in Ferguson) or quadrangles (as in the present discussion of the Fels case). The analysis of the construction of Fels sculptures for quadrangles as discussed above may be extended to cover $n$-sided curvilinear polygons (whether they are planar or not) and a finite sequence, $t_k^*$, of special parameter values, $t^*$,  in the interval $[0,a]$. 
\begin{figure}[ht!]
\centering
\includegraphics[width=88.5mm]{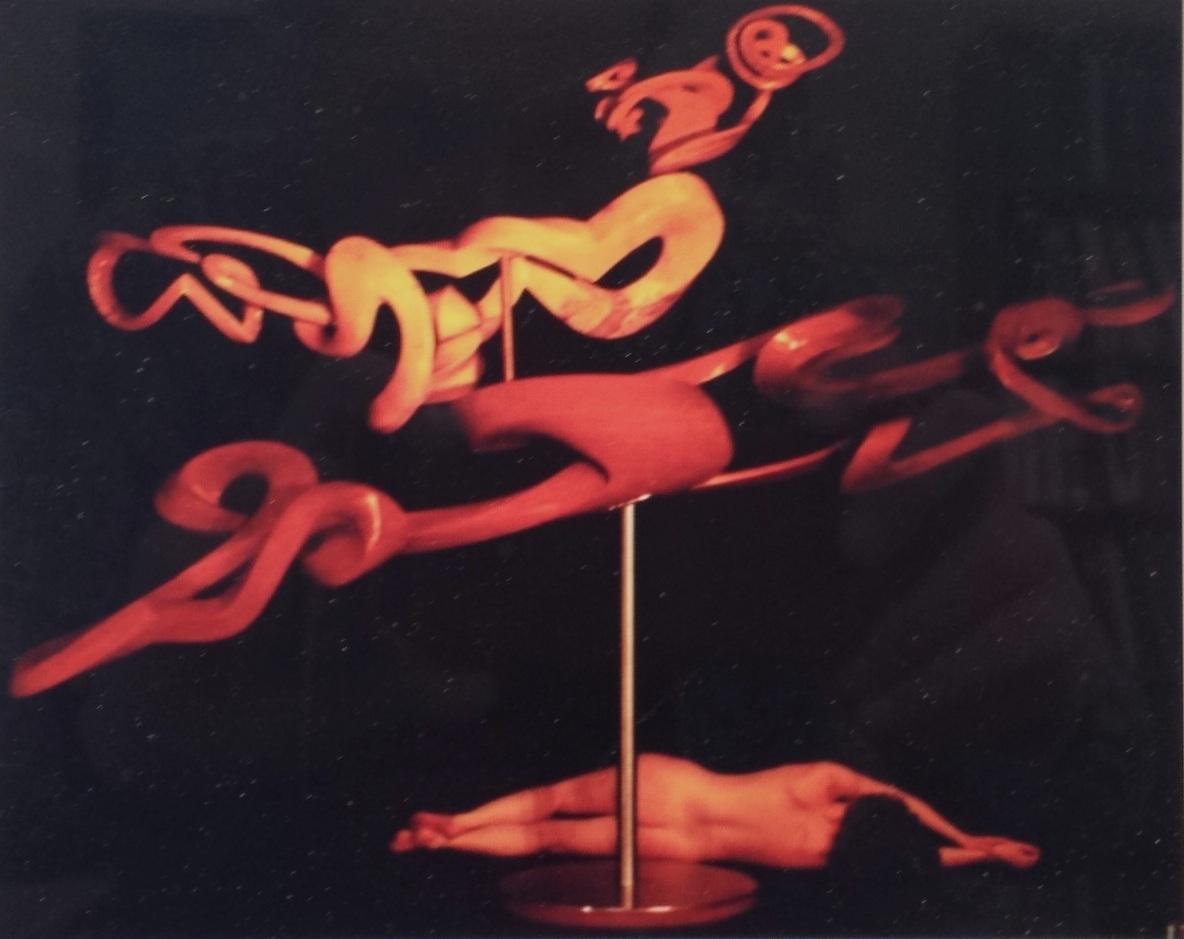}
\includegraphics[width=50mm]{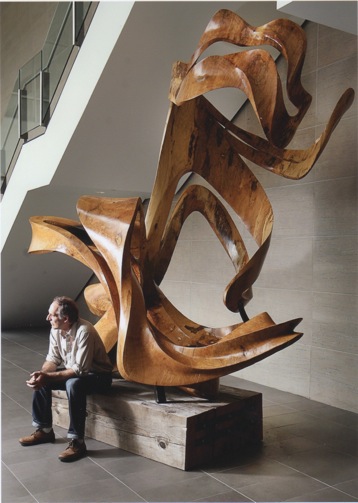}
\caption{Complex Fels sculptures}
\label{Fels16}
\end{figure}
The same rules apply since all that needs to be done is replace the word {\it quadrangle} by an arbitrary $n$-sided, possibly non-planar, polygon and then replace the one value of $t^*$ by a finite sequence of such between each of which there is a special (initially fixed) clockwise or counterclockwise rotation of the polygons during the continuous deformation. In addition, for a given Fels sculpture there is always the requirement that there is a unique cross-section having minimal area and a unique cross-section having maximal area. Of course, in this case, there must exist cross-sections of the sculpture having a locally minimal (resp. locally maximal) area.

\begin{figure}[ht!]
\centering
\includegraphics[width=51mm]{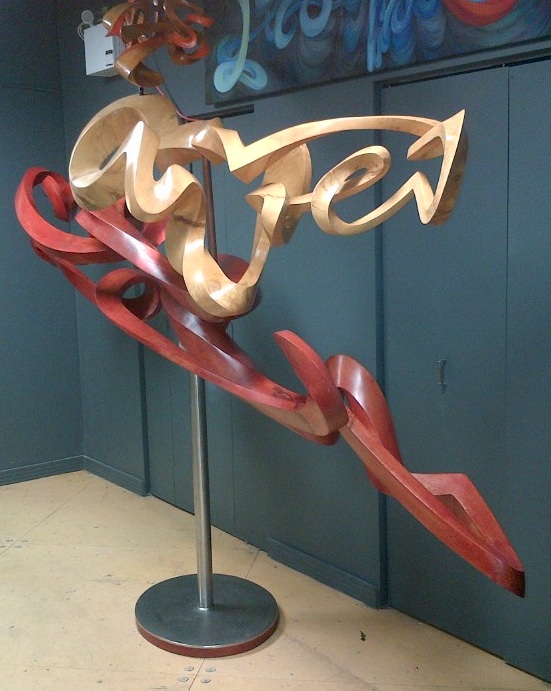}
\includegraphics[width=46mm]{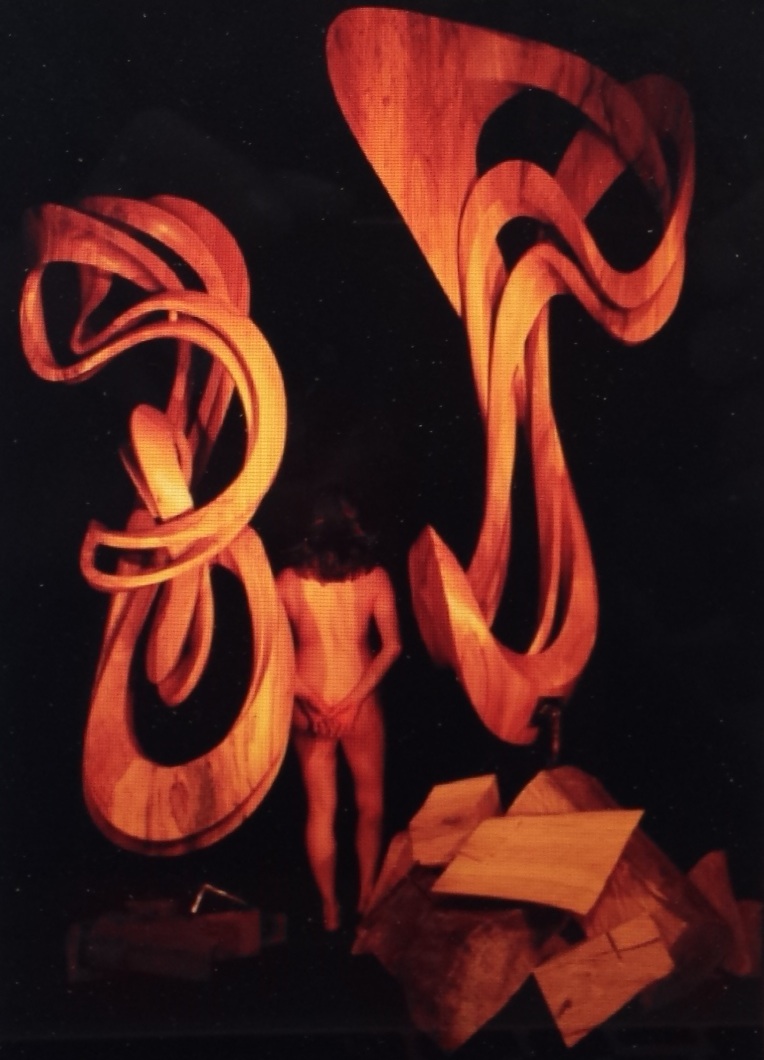}
\includegraphics[width=49mm]{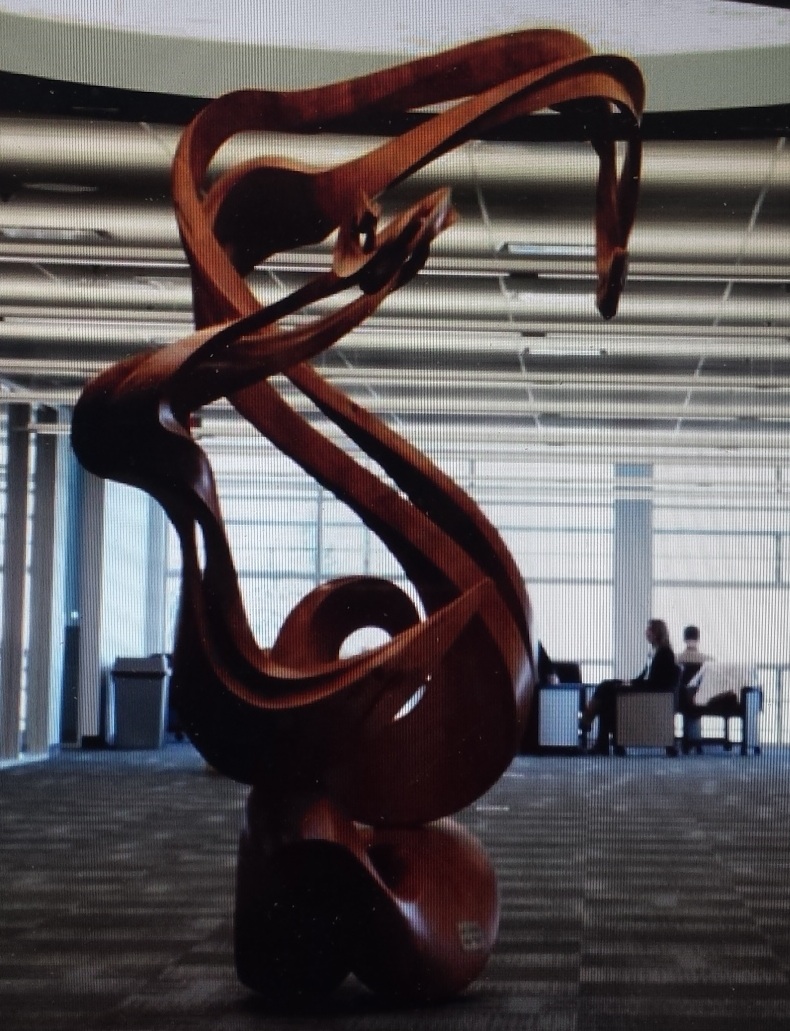}
\caption{Variations on the theme of Fels sculptures using a sequence of times and tangential groupings of simpler Fels sculptures}
\label{Fels15}
\end{figure}

In another development, Fels considered the possibility of the center of mass curve intersecting itself {\it tangentially} thus forcing the quadrangular structure of the sculpture to create solids that are no longer toroidal and of essentially arbitrary genus (Figures~\ref{Fels16},~\ref{Fels15}).

Examples of curvilinear quadrangles in the setting of a {\it sequence} of special parameter values, $t_k^*$, between each of which the rotations vary from clockwise to counterclockwise by a specific angle, are shown below in Figure~\ref{Fels15}.

\section*{Virtual Sculpting}
Select two unique quadrangles and use the respective centroids as nodes wirelessly connected to two gloves for your hands. As you move your hands, you shape the sculpture in real time while the language, as controlled by the software, modifies the sculptural form accordingly. This can be as complex or as simple as you like. This process enables ongoing exploration and changes of expression of the form through human gesture. For example, interfacing with music is definitely possible generating a projected three-dimensional ‘light organ’. The human spirit, therefore, can be personified through virtual sculpting. 

The capture of human movement with wearable sensors is well suited to creating real time virtual sculpting using the language mathematically described in this paper. We propose a simple example as follows:

\begin{figure}[ht!]
\centering
\includegraphics[width=73.5mm]{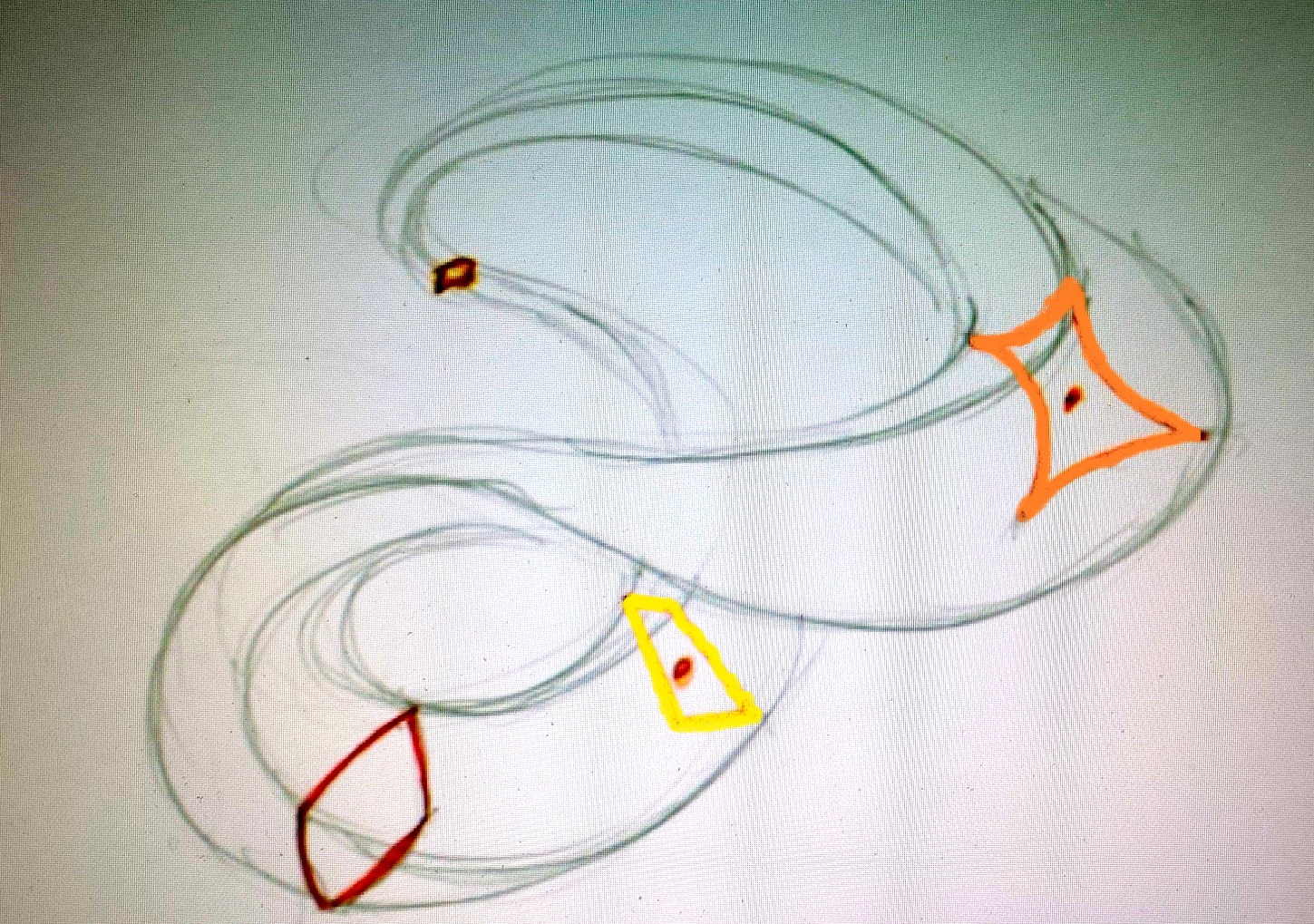}
\includegraphics[width=70mm]{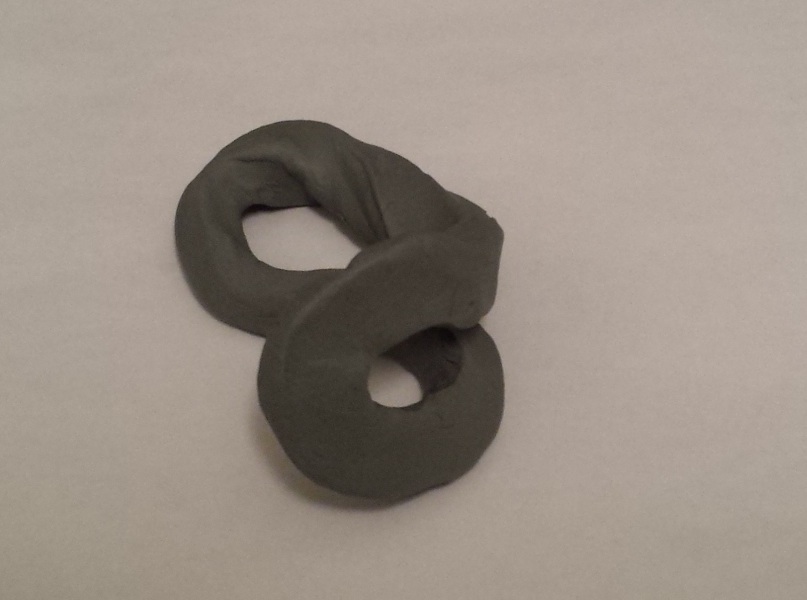}
\caption{An example of virtual sculpting and a clay mock-up}
\label{Fels14}
\end{figure}

In Figure~\ref{Fels14} default red and orange quadrilaterals are preprogrammed as stationary in time and space relative to each other as projected on a screen. There are two gloves wherein: Each glove has a sensor in the palm representing the centroid of an orange quadrilateral. Each glove has 4 sensors, one on the thumb and three on its nearest fingertips, representing the vertices of the respective orange quadrilateral. The person modifies the orange quadrilaterals though virtual space and real time through the following rules:

Moving the hands in space changes the locations of the orange quadrilaterals relative to each other and to the default red and yellow quadrilaterals. These variations are programmed to follow the sensor in the palm of each glove. The vertices of the orange quadrilaterals are programmed to be expanded away or towards their respective centroids by opening and closing the fingers. Note the language restricts the area of each of the orange quadrilaterals to larger than the yellow quadrilateral. The programming based on the language would fill in all other quadrilaterals between the default red, yellow, and orange quadrilaterals, in the smoothest curve, as per the language of form defined here.

In this simple example, moving your hands shapes the sculpture in real time because the language, as governed by the mathematics, modifies the sculptural form accordingly. This can be as complex or as simple as you like. This process enables ongoing exploration and changes of expression of form through human gesture.

\section*{Sculptures in hyperspace}

Finally the restriction to our own visually motivated three-dimension
al world is for esthetic purposes only. There is no real need to make such a restriction on mathematical grounds since the construction there may be formulated for any number of finite or even infinitely many dimensions. To bring these sculptures back {\it down} to our world all that is needed is a suitable three-dimensional projection. Thus, it is conceivable that a Fels-like-sculpture in four dimensional space may generate many, many different projections each being a three-dimensional Fels-sculpture in its own right. In this case, the process begins with an arbitrary curve in four-dimensional space (containing the centers of mass of our four, or even lower dimensional,  dimensional quadrangles) and an arbitrary $n$-sided polygon where $n$, the number of sides, is fixed. This can  also be called an $n$-sided {\it $m$-polytope}, where $m$ is the dimension in which it resides.  Using analytic geometry one can write down the equations of a typical continuous deformation of such a polytope and create mathematically real Fels sculptures in four dimensional space. Of course, we cannot visualize such a sculpture but we {\it can} see its projections back into our space. 

{\bf EXAMPLE}

\begin{figure}[ht!]
\centering
\includegraphics[width=70mm]{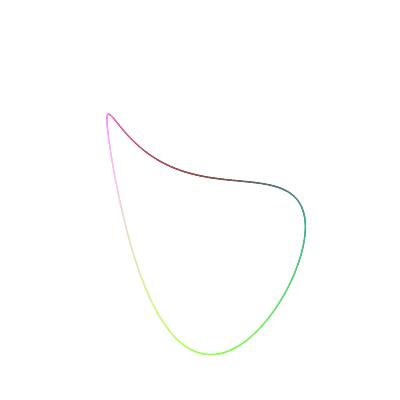}
\includegraphics[width=70mm]{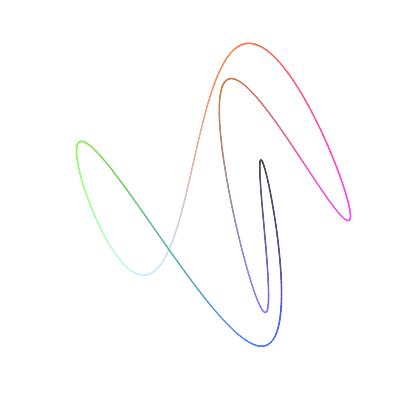}
\includegraphics[width=70mm]{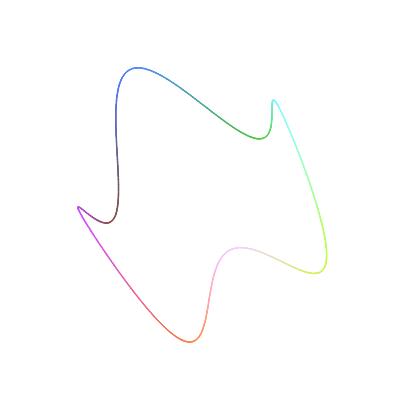}
\includegraphics[width=70mm]{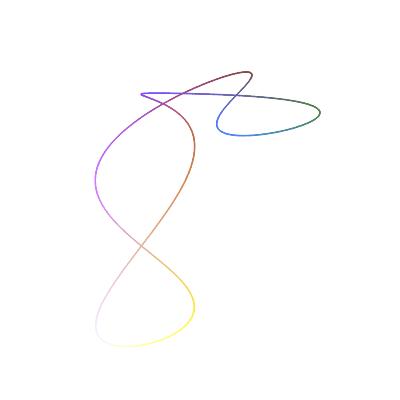}
\caption{Wireframes of four three dimensional projections (on a plane) of a four-dimensional Fels sculpture. Clockwise from the top left we see the $xyz$, $yzw$, $zwx$, and $wxy$ projections of the four dimensional curve of this example on a plane.}
\label{Fels17}
\end{figure}

We consider the closed four-dimensional curve given parametrically by
$$\mathbf{r}(t) = (\cos t, \sin t, \cos 2t +\cos t, \sin 4t), \quad 0 \leq t \leq 2\pi.$$
Of course, this curve cannot be sketched but its various {\it projections} in three dimcnsional (humanly visible) space can be manifested easily using software suited for this purpose (e.g., Maple\textsuperscript{\textregistered}). This four dimensional (center of mass) curve may be used to create the wireframe of many different three dimensional Fels sculptures (see Figure~\ref{Fels17}) using the ideas presented earlier. Here we show its four main projections on the three-dimensional $xyz$, $yzw$, $zwx$ and $wxy$ spaces, and then those projections onto the plane of this reading surface. Observe that the four curves in Figure~\ref{Fels17} are not self-intersecting when viewed in three dimensions, though they are when viewed in two (because of the nature of projections). Once the wireframes are constructed the rest is a matter of artistic license subject to the basic principles outlines in this work.

\end{document}